\documentclass[a4paper, 11pt]{article} 
\usepackage{amsfonts,amsmath,amssymb,amscd,color,stmaryrd,latexsym}
\usepackage[all]{xy} 
\usepackage[T1]{fontenc}

\def\ker{\mathop{\rm ker}\nolimits}

\def\coker{\mathop{\rm coker}\nolimits}

\def\Re{\mathop{\rm Re}\nolimits}
\def\Im{\mathop{\rm Im}\nolimits}

\def\mod{\mathop{\rm mod}\nolimits}
\def\note#1{\marginpar{\raggedright\if@twoside\ifodd\c@page\raggedleft\fi\fi\sf\scriptsize RMK: #1}}

\newcommand{\ind}{\mathrm{index}\,}

\newcommand{\CP}{{\mathbb C}{\mathbb P}}

\newcommand{\mf}{\mathfrak}
\newcommand{\mb}{\mathbf}
\newcommand{\mc}{\mathcal}

\newcommand{\R}{\mathbb{R}}
\newcommand{\Z}{\mathbb{Z}}
\newcommand{\C}{\mathbb{C}}

\newcommand{\Oc}{\mathbb{O}}
\newcommand{\Hh}{\mathbb{H}}

\newcommand{\wrt}{with respect to }

\newcommand{\loc}{\stackrel{\mathrm{loc}}{=}}

\newenvironment{rmk}{\begin{trivlist}\item[]{\bf Remark:}\setlength{\parindent}{0pt}}{\end{trivlist}}
\newenvironment{ex}{\begin{trivlist}\item[]{\bf Example:}\setlength{\parindent}{0pt}}{\end{trivlist}}
\newenvironment{prf}{\begin{trivlist}\item[]{\bf Proof: }}{\hfill $\blacksquare$ \end{trivlist}}

\newtheorem{thm}{Theorem}[section]
\newtheorem{definition}[thm]{Definition}
\newtheorem{prp}[thm]{Proposition}
\newtheorem{lem}[thm]{Lemma}
\newtheorem{cor}[thm]{Corollary}

\begin{document}
\title{Deformations of associative submanifolds with boundary}
\author{Damien Gayet and Frederik Witt}
\date{}
\maketitle

\centerline{\textbf{Abstract}}
Let $M$ be a topological $G_2$--manifold. We prove that the space of infinitesimal associative deformations of a compact associative submanifold $Y$ with boundary in a coassociative submanifold $X$ is the solution space of an elliptic problem. For a connected boundary $\partial Y$ of genus $g$, the index is given by $\int_{\partial Y}c_1(\nu_X)+1-g$, where $\nu_X$ denotes the orthogonal complement of $T\partial Y$ in $TX_{|\partial Y}$ and $c_1(\nu_X)$ the first Chern class of $\nu_X$ with respect to its natural complex structure. Further, we exhibit explicit examples of non--trivial index.
 
\bigskip


\smallskip

\section{Introduction}
The deep and rich interplay between geometry and algebra on manifolds with $G_2$--structure is reflected in the existence of special submanifolds, namely {\em associative} ones of dimension~$3$ and {\em coassociative} ones of dimension~$4$. These are particular instances of Harvey's and Lawson's {\em calibrated submanifolds}~\cite{hala82}, a notion which also embraces complex submanifolds of a K\"ahler manifold or special lagrangian submanifolds of a Calabi--Yau manifold. McLean~\cite{mc98} proved that the infinitesimal coassociative deformations of a coassociative $X$ is an unobstructed elliptic problem. The dimension of the moduli space is $b^2_{+}(X)$, i.e.\ the dimension of self--dual harmonic 2--forms on $X$. For associative submanifolds, the problem, though still elliptic, is more involved. Firstly, existence of smooth deformations is, as in the case of complex submanifolds, in general obstructed. Secondly, the {\em virtual} dimension, that is, the index of the elliptic equation, is always zero on dimensional grounds, so that no prediction on the existence of infinitesimal deformations can be made. This result was extended to arbitrary manifolds with topological $G_2$--structure (whose holonomy is not necessarily contained in $G_2$) by Akbulut and Salur~\cite{aksa07}, who also address smoothness and compactness issues of the deformation spaces.

\bigskip

On the other hand, on symplectic manifolds one is naturally led to study the moduli space of (pseudo--)holomorphic curves with boundary in a lagrangian submanifold~\cite{fl88},~\cite{gr85}. In physics, Aganagic and Vafa translated this boundary problem for special lagrangians of a Calabi--Yau into an open string problem~\cite{agva00}, following Witten's use of the moduli space of complex curves in the stringy world~\cite{wi91}. Now taking a Calabi--Yau $3$--fold $K$ times a circle yields a natural riemannian manifold $M=K\times S^1$ with holonomy contained in $G_2$. Moreover, holomorphic curves and special lagrangians times a circle give examples of associative and coassociative submanifolds in $M$. In this way, the duality of complex versus special lagrangian submanifolds is matched by the duality of associative versus coassociative submanifolds in a holonomy $G_2$--manifold. On one side, this hints at the existence of a $G_2$--analogue of Floer theory as conjectured in~\cite{dose09}. Further, it makes (co--)associatives play a key r\^ole in attempts to set up topological $M$--theory~\cite{btz06}, \cite{dggv05}. It is therefore natural to study deformations of (co--)associatives with boundary. Inspired by the work of Butscher~\cite{bu01}, who investigated deformations of special lagrangians with boundary on a symplectic, codimension $2$ submanifold inside some compact Calabi--Yau, Kovalev and Lotay investigated in a recent paper the analogous problem for manifolds with closed $G_2$--structures, where a compact coassociative has its boundary in a fixed, codimension $1$ submanifold~\cite{kolo07}. 

\bigskip 

In this paper, we consider an associative $Y$ inside topological $G_2$--manifolds and study the space of infinitesimal associative deformations of $Y$ with boundary inside a fixed coassociative $X$. We identify this space with the solutions of an elliptic boundary value problem whose index is given for connected boundary $\partial Y$ of genus $g$ by
$$
\ind(X,Y)=\int_{\partial Y}c_1(\nu_X)+1-g.
$$ 
Here, $\nu_X$ denotes the orthogonal complement of $T\partial Y$ in $TX_{|\partial Y}$ and $c_1(\nu_X)$ the first Chern class of $\nu_X$ with respect to a natural complex structure we are going to define below. Further, we extend this result to $4$--dimensional submanifolds $X$ which do not contain any associative. In a sense, this class of submanifolds inside a (topological) $G_2$--manifolds forms the natural counterpart of totally real submanifolds inside (almost) complex manifolds. Finally, assuming that $Y$ is an embedded $3$--disk, we associate with $Y$ an element $\mu_{G_2}(\partial Y)\in\pi_2\big(G_2/SO(4)\big)\cong\Z_2$, which is best thought of as a $G_2$--analogue of the Maslov index. Under suitable identifications, we show that
$$
\mu_{G_2}(X,Y)=\int_{\partial Y}c_1(\nu_X)\mod2=\big(\ind(X,Y)+1\big)\mod2.
$$
Explicit examples of pairs $(X,Y)$ with non--trivial index will be given in Section~\ref{indexa}. In particular, we shall construct compact pairs $(X,Y)$ inside compact holonomy $G_2$--manifolds using Joyce's method~\cite{jo00}. 

\bigskip

The techniques we use are the standard ones from PDE theory; our reference is~\cite{bowo93} whose conventions we shall follow throughout this paper.
%
%
%
%
%
\section{The group $G_2$}
\label{groupg2}
We start by recalling some classical facts about $G_2$ (cf.\ for instance~\cite{aksa07},~\cite{br06},~\cite{hala82} and~\cite{jo00}).
%
%
\subsection{Linear algebra}
\label{octonions}
The octonions define a real $8$--dimensional, non--associative division algebra $\Oc=\Hh\oplus e\Hh$ generated by $\langle\mb{1},i,j,k,e,e\cdot i,e\cdot j,e\cdot k\rangle$. Taking these generators as an orthonormal basis induces an inner product $\langle\cdot\,,\cdot\rangle$ on $\Oc$ compatible with the algebra structure. Further, we obtain a {\em vector cross product} taking values in the imaginary octonions $\Im\Oc=\langle\mb{1}\rangle^{\perp}\cong\R^7$ defined by
$$
u\times v=\Im(\overline{v}\cdot u).
$$
Here, $\overline{v}$ is the natural conjugation which sends $v\in\Im\Oc$ to $-v$. The term cross product is justified by the properties $u\times v=-v\times u$ and $|u\times v|=|u\wedge v|$.
Over $\R^7$, this yields the $3$--form
$$
\varphi_0(u,v,w)=\langle u\times v,w\rangle,
$$
which \wrt the orthonormal basis $e_1=i,\,e_2=k,\ldots,e_7=e\cdot k$ can be written as\footnote{This is the convention adopted in~\cite{aksa07},~\cite{br06} and~\cite{jo00}.}
\begin{equation}\label{standform}
\varphi_0=e^{123}+e^1\wedge(e^{45}+e^{67})+e^2\wedge(e^{46}-e^{57})+e^3\wedge(-e^{47}-e^{56}).
\end{equation}

\bigskip

We refer to any basis $\{e_j\}$ such that $\varphi_0$ is of the form~(\ref{standform}) as a $G_2$--{\em frame}, since the stabiliser of $\varphi_0$ inside $GL(7)$ is the real algebraic Lie group $G_2$, which is of dimension $14$. Conversely, any $G_2$--invariant form $\varphi\in\Lambda^3\R^{7*}$ induces a 
positive definite inner product $\langle\cdot\,,\cdot\rangle_{\varphi}$ and a 
cross product $\times_{\varphi}$ as follows. Firstly, with $\varphi$ we can associate 
a volume form $\mu_{\varphi}$ (which is somehow difficult to write down explicitly, cf.\ the appendix in~\cite{hi01}). Then we define
\begin{equation}\label{metriccross}
\langle u,v\rangle_{\varphi}=\big((u\llcorner\varphi)
\wedge(v\llcorner\varphi)\wedge\varphi\big)/6\mu_{\varphi},\quad
\langle u\times_{\varphi} v,w\rangle_{\varphi}=\varphi(u,v,w).
\end{equation}

\bigskip

Next, we consider the {\em associator}
$$
[u,v,w]=\frac{1}{2}\big((u\cdot v)\cdot w-u\cdot(v\cdot w)\big)
$$
which is totally skew--symmetric. The associated $\Im\Oc$--valued $3$--form in $\Lambda^3\R^{7*}\otimes\R^7$ will be written $\chi_0$. For this form, we have the important identity
\begin{equation}\label{chirule}
\chi_0(u,v,w)=-u\times(v\times w)-\langle u,v\rangle w+\langle u,w\rangle v.
\end{equation}
In particular, we find $u\times(u\times a)=-|u|^2\cdot a$ if $a$ is orthogonal to $u$. 

\bigskip

Finally we define a $4$--form over $\R^7$ by
$$
\psi_0(u,v,w,x)=\frac{1}{2}\langle u,[v,w,x]\rangle.
$$
This form actually coincides with the Hodge dual of $\varphi_0$, so that in a $G_2$--frame $\{e_j\}$, 
$$
\psi_0=\star\varphi_0=-e^{12}\wedge(e^{47}+e^{56})-e^{13}
\wedge(e^{46}-e^{57})+e^{23}\wedge(e^{45}+e^{67})+e^{4567}.
$$
%
%
\subsection{$G_2$--manifolds}
\label{g2manif}
Next consider a $7$--dimensional manifold $M$. If the structure group $GL(7)$ reduces to $G_2$, we say that $M$ carries a $G_2$--{\em structure} and refer to $M$ as a {\em topological $G_2$--manifold} or simply as a $G_2$--{\em manifold} for sake of brevity. In this case, there exists a $3$--form $\varphi$ on $M$ such that the associated $G_2$--principal frame bundle consists of isomorphisms between $(T_xM,\varphi_x)$ and $(\R^7,\varphi_0)$ for $x\in M$. By an {\em abus de langage}, we refer to the defining $3$--form $\varphi$ itself as the $G_2$--structure. Then there exists a vector cross product $\times=\times_{\varphi}$ on $TM$, inducing the structure of $\Im\Oc$ on any tangent space $T_xM$. Moreover, formulae~(\ref{metriccross}) give rise to a globally defined riemannian metric $g=g_{\varphi}$ with Levi--Civita connection $\nabla^g$. Similarly, there are global counterparts $\psi=\star_{\varphi}\varphi\in\Omega^4(M)$ and $\chi\in C^{\infty}(M,\Lambda^3T^*M\otimes TM)$ of $\psi_0$ and $\chi_0$ respectively. A $G_2$--manifold is said to be {\em torsion--free} if $\nabla^g\varphi=0$. By a result of Fernandez and Gray~\cite{fegr82}, torsion--freeness is equivalent to $d\varphi=d\psi=0$. Note that in this case, the holonomy of $g$ is contained in $G_2$. This is tantamount to saying that there exist coordinates around each $x\in M$ such that $\varphi(x)=\varphi_0+O(|x|^2)$. In the sequel, we say that a torsion--free $G_2$--manifold is a {\em holonomy $G_2$--manifold}, if the holonomy actually equals\footnote{Note that some authors do not make this distinction and refer to any torsion--free $G_2$--manifold as a holonomy $G_2$--manifold or even simply as a $G_2$--manifold.} $G_2$. 

\bigskip

An oriented $3$--dimensional submanifold $Y$ is called {\em associative} if the pull--back of $\varphi$ to $Y$ is equal to the induced riemannian volume form. Equivalently, the pull--back of $\chi$ to $Y$ is identically zero. Furthermore, an oriented $4$--dimensional submanifold $X$ is called {\em coassociative} if the pull--back of $\psi$ to $X$ is equal to the induced riemannian volume form. Equivalently, the pull--back of $\varphi$ to $X$ is identically zero. As shown by Harvey and Lawson, (co--)associative manifolds have the important property of being homologically volume minimising if the form $\varphi$ (respectively $\psi$) is closed. In particular, following their language, $\varphi$ and $\psi$ define {\em calibrations}, and (co--)associative submanifolds are {\em calibrated}.
%
%
%
%
%
\section{The geometry of the deformation problem}
Let $M$ be a $G_2$--manifold and $Y\subset M$ a compact associative whose boundary $\partial Y$ is contained in a fixed $4$--submanifold $X\subset M$. 
We wish to describe the space of infinitesimal deformations of $Y$ in the class of associatives with boundary in $X$, that is, the {\em Zariski tangent space} $T^{\mathrm{Zar}}\mf{M}_{X,Y}$ of
$$
\mc{M}_{X,Y}=\{Y'\,|\,Y'\mbox{ associative isotopic to }Y\mbox{ with }\partial Y'\subset X\}.
$$
%
%
\subsection{The closed case}
\label{closedcase}
Our starting point is the closed case $\partial Y = \emptyset$, that is we first discuss the Zariski tangent space $T^{\mathrm{Zar}}\mf{M}_Y$ of $\mc{M}_Y=\{Y'\,|\,Y'\mbox{ associative isotopic to }Y\mbox{ with }\partial Y'=\emptyset\}$. In the case of associatives inside {\em torsion--free} $G_2$--manifolds, this problem was settled by McLean~\cite{mc98}. Later on, his result was generalised by Akbulut and Salur~\cite{aksa04a},~\cite{aksa07} to arbitrary $G_2$--manifolds. We outline these results following the approach of~\cite{aksa07}. 

\bigskip

Let $(M,\varphi,\nabla)$ be a $G_2$--manifold endowed with a compatible connection $\nabla$, i.e.\ $\nabla\varphi\equiv0$. Such connections always exist (cf.\ for instance~\cite{br06}), but are not necessarily unique. In particular, $\nabla$ preserves the metric and therefore induces a connection on the normal bundle $\nu\to Y$, which we also denote by $\nabla$. Since $Y$ is compact, any deformation $Y\to Y_t$, $t\in(-\epsilon,\epsilon)$ can be assumed to be normal: Reparametrising by a time--dependent diffeomorphism if necessary, one can achieve that $\partial Y_t(p)/\partial t \in\nu_p$. Consequently, nearby (i.e.\ $C^1$--close) $Y_t$ can be identified with sections $s_t\in C^{\infty}(Y,\nu)$ via the exponential map. The associative $Y$ then corresponds to $s_0=$ zero--section of $\nu$. Pulling $TY_t$ back to $TM_{|Y}$ through parallel transport \wrt $\nabla$ along the curves $\lambda\mapsto\exp_p\big(\lambda s_t(p)\big)$, $p\in Y$, we obtain a curve in the fibre $G_3(TM)_p=G_3(T_pM)$. For $t=0$, it passes through $E_p:=T_pY$, which is an element in $G_{\varphi}(TM)_p$, the subset of calibrated $3$--planes in $G_3(TM)_p$. The derivative of this curve at $0$ can be identified with $\nabla\dot{s}_0(p)$, and the linearised condition is thus $\nabla\dot{s}_0(p)\in T_{E_p}G_{\varphi}(M)\subset T_{E_p}G_3(M)\cong E_p^*\otimes\nu_p$. Moreover, the vector cross product gives rise to an exact sequence
$$
0\to T_{E_p}G_{\varphi}(TM)\to T_{E_p}G_3(TM)\cong E_p^*\otimes\nu_p\stackrel{\times}{\to}\nu_p\to0.
$$
Picking an orthonormal basis $e_1,\,e_2,\,e_3$ of $E_p$ we obtain the condition $\nabla\dot{s}_0(p)\in T_{E_p}G^{\varphi}_3(M)$ if and only if $\sum e_i\times\nabla_{e_i}\dot{s}_0(p)=0$. An invariant formulation of this equation can be given in terms of a {\em Dirac operator}\footnote{For sake of brevity, we refer to any operator $\mb{D}$  as a Dirac operator if it is of {\em Dirac type}, that is, the principal symbol of $\mb{D}^2$ satisfies $\sigma(\mb{D}^2)(p,\xi)=||\xi||^2$, cf.\ also~\cite{bowo93}, Section 3.}. The fibre bundles $TY$ and $\nu$ are associated  with an $SO(4)$--representation. Further, since $Y$ is spinnable the principal $SO(4)$--bundle can be lifted to a principal $Spin(4)$--bundle. Then $\nu$ is associated with the tensor product of a spin representation and some other representation of $Spin(4)$. As a result we may regard $\nu$ as a {\em twisted spinor bundle}. Under this identification, the operator $\sum e_i\times\nabla$ becomes the Dirac operator of $(\nu,\nabla)$. Summing up, we arrive at a generalised version of McLean's

\begin{thm}\hspace{-4pt}{\rm(\cite{aksa07},~\cite{mc98})}\label{mclean}
Let $(M,\varphi)$ be a $G_2$--manifold and $Y\subset M$ an associative. Then $T^{\mathrm{Zar}}\mf{M}_Y$ can be identified with the kernel of a twisted Dirac operator $\mb{D}^{\nabla}:C^{\infty}(Y,\nu)\to C^{\infty}(Y,\nu)$ taken with respect to a connection $\nabla$ induced by a compatible connection of $(M,\varphi)$. In particular, we obtain the natural Dirac operator on $\nu$ \wrt $\nabla^g$ if $(M,\varphi)$ is torsion--free.
\end{thm}

In the sequel, we denote the Dirac operator $\mb{D}^{\nabla}$ of $Y$ simply by $\mb{D}$.

\begin{rmk}
The {\em Dirac equation} $\mb{D}s=0$ is elliptic, and as a consequence, the {\em virtual dimension} of $\mf{M}_Y$, that is, $\dim T^{\mathrm{Zar}}\mf{M}_Y$, is finite. Furthermore, $Y$ being odd--dimensional, we have $\ind(\mb{D})=\dim\ker(\mb{D})-\dim\coker(\mb{D})=0$. In generic situations where one expects the cokernel to vanish, we would get $\dim T^{\mathrm{Zar}}\mf{M}_Y=0$, and associatives would be rigid objects. However, we have no a priori control on the virtual dimension in terms of topological datum. This stands in sharp contrast to the deformations of a coassociative $X$, where McLean showed that the (actual) dimension of the (smooth) moduli space $\mf{M}_X$ is $b^2_+(X)$, the dimension of harmonic self--dual $2$--forms~\cite{mc98}.
\end{rmk}
%
%
\subsection{The geometry on the boundary}
\label{boundgeo}
Next assume that the associative $Y$ has a non--empty boundary $\partial Y$ inside a fixed coassociative $X$. We first need to understand the geometry on the boundary of $Y$.

\bigskip

Fix a collar neighbourhood $\mc{C}\cong\partial Y\times[0,\epsilon)$ of $\partial Y$ and let $u$ denote the inward pointing unit vector field defined on $\mc{C}$. As before, $\nu\to Y$ denotes the normal bundle, as well as its restriction to $\partial Y$. As remarked in~\cite{aksa04a} and~\cite{aksa07}, $\nu_{|\mc{C}}$ carries a hermitian structure near the boundary induced by $u$, namely
$$
J:\nu\to\nu,\quad J(x)=u\times x.
$$
This acts indeed as an isometry with respect to $g$, as
$$
g(Ja,Jb)=\varphi(u,a,u\times b)=-g\big(u\times(u\times b),a\big)=g(a,b)
$$
for any $a,\,b\in\nu_{|\mc{C}}$. Let $\nu_X\subset TX_{|\partial Y}$ denote the orthogonal complement of $T\partial Y $ in $TX_{|\partial Y}$, i.e.
$$
TX_{|\partial Y}=T\partial Y\oplus\nu_X.
$$

\begin{lem}\label{boundary} 
For the bundle $\nu\to\partial Y$ the following holds: 
\begin{enumerate}
\item The bundle $\nu_{X}$ is contained in $\nu$ and is stable under $J$.
\item The orthogonal complement $\mu_X$ of $\nu_X$ in $\nu$ is also stable under $J$.
\item Viewing  $T\partial Y$, $\nu_X$ and $\mu_X$ as $J$--complex bundles, we have 
$$
\overline{\mu}_X\cong\nu_X\otimes_{\C}T\partial Y,
$$
that is $\mu_X^{0,1}\cong\nu_X^{1,0}\otimes T^{1,0}\partial Y\cong\nu^{1,0}_X\otimes\overline{K}_{\partial Y}$, where $K_{\partial Y}$ is the canonical line bundle over $\partial Y$.
\end{enumerate}
\end{lem}
\begin{prf}
Let us fix a local orthonormal basis $(u,v,w)$ of $TY_{|\partial Y}$ by choosing a local unit vector field $v\in T\partial Y$, i.e.
$$
TY_{|\partial Y}\loc\langle u,v,w\rangle.
$$
The vector field $w$ is defined by $w=u\times v$, which lies in $T\partial Y$ in virtue of the associativity of $Y$. If $a\in\nu_X$, then $g(a,u)=0$, for $v\times w=u$ and $\varphi(v,w,a)=0$, $X$ being coassociative. Clearly, the vectors $a\times v$ and $a\times w$ are orthogonal to $v$ and $w$ as well as to $u$, since 
$$
g(a\times v,u)=\varphi(a,v,u)=-g(u\times v,a)=-g(w,a)=0,
$$
and similarly for $a\times w$. Hence $a\times v$, $a\times w\in\nu$. Further, these vectors are orthogonal to $TX$, for $a,\,v,\,w\in TX$ and $X$ is coassociative, so that for $n\in\nu_X$ we find $g(a\times v,n)=\varphi(a,v,n)=0$ etc.. Hence
$$
\mu_X\loc\langle a\times v,a\times w\rangle.
$$
As a consequence, $u\times a\in\nu$ is orthogonal to $\mu_X$ \big(for $g(u\times a, a\times v)=\varphi(u,a,a\times v)$ etc.\big), so that
$$
\nu_X\loc\langle a,Ja=u\times a\rangle,\quad TX_{|\partial Y}\loc\langle v,w,a Ja\rangle.
$$
This shows that $\nu_X$ is stable under $J$. On the other hand, $g\big(u\times (a\times v),a\big)=\varphi\big(u,(a\times v),a\big)=0$ and similarly $g\big(u\times(a\times v),u\times a\big)=0$, hence $u\times (a\times v)\in\mu_X$ which shows that $\mu_X$ is also stable under $J$. 

\bigskip

The Riemann surface structure on $\partial Y$ is induced by the hermitian structure $J=u\times$, for $Jv=u\times v=w$ and $Jw=u\times w=-v$ (to keep notation tight we abuse notation and also write $J$ for the endomorphism on $T\partial Y$ induced by $u\times$). The map
$$
a\otimes y\in \nu_X\otimes_{\C} T\partial Y\mapsto a\times y\in\overline{\mu}_X,
$$
where we now view $\nu_X$ and $\overline{\mu}_X$ as complex line bundles 
via $J$, is well--defined and a {\em real} bundle isomorphism. It remains to see that it is complex--linear, i.e.
$$
Ja\times y=a\times Jy=-J(a\times y).
$$
This is equivalent to $(u\times a)\times y$ and $a\times(u\times y)$ 
being equal to $-u\times(a\times y)$. But this follows from~(\ref{chirule}) and the skew--symmetry of $\chi$.
\end{prf}

Let $\mb{B}:C^{\infty}(\partial Y,\nu)\to C^{\infty}(\partial Y,\mu_X)\subset C^{\infty}(\partial Y,\nu)$ be the orthogonal projector\footnote{By an {\em orthogonal projector} we understand an operator $\mb{P}$ of order $0$ satisfying $\mb{P}=\mb{P}^2=\mb{P}^*$.} taking smooth sections of $\nu=\nu_X\oplus\mu_X$ to $\mu_X$. As a corollary to the generalised version of McLean's theorem as given in Theorem~\ref{mclean} and Lemma~\ref{boundary}, we obtain:

\begin{cor}\label{defequs}
The Zariski tangent space of $\mc{M}_{X,Y}$ can be identified with solutions of the system
\begin{equation}\label{defprobreal} 
\mb{D}s=0\,,\quad\mb{B}(s_{|\partial Y}) =0,\quad s\in C^{\infty}(Y,\nu).
\end{equation} 
\end{cor}
%
%
%
%
%
\section{Ellipticity and index}
In view of applying the standard machinery of index theory for manifolds with boundary, we consider the complexification of equation~(\ref{defprobreal}), namely 
\begin{equation}\label{defprob}
(\mb{D}^{\C}\oplus\mb{B}^{\C})(\sigma\oplus\sigma_{|\partial Y})=0,\quad s\in C^{\infty}(Y,\nu^{\C})
\end{equation} 
with operators $\mb{D}^{\C}$ and $\mb{B}^{\C}$ extended to the complexified bundles $\nu^{\C}=\nu\otimes\C$ and $\mu_X^{\C}=\mu_X\otimes\C$. However, we shall not distinguish between the original operators and their complexification in the sequel for sake of keeping notation tight. To show that the kernel is finite--dimensional requires a suitable ellipticity condition. We first introduce the {\em Calder\'on projector} $\mb{Q}_{\mb{D}}$ associated with a Dirac operator $\mb{D}$ (cf.~\cite{bowo93}, Thm. 12.4). This is a zero order pseudo--differential operator 
$$
\mb{Q}_{\mb{D}}:C^{\infty}(\partial Y,\nu)\to\mc{H}(\mb{D})=\{s_{|\partial Y}\,|\,s\in C^{\infty}(Y,\nu),\,\mb{D}s=0\}\subset C^{\infty}(\partial Y,\nu)
$$ 
mapping smooth sections of $\nu$ over $\partial Y$ to the space of Cauchy data\footnote{We are glossing over some technical details such as the passing to the ``closed double'' $M=Y\cup_{\partial Y}Y$, cf.\ Chapters 9, 11 and 12 in~\cite{bowo93}.} of $\mb{D}$.

\begin{definition}\label{locbond}\hspace{-4pt}{\rm (cf.~\cite{bowo93}, Def.\ 18.1)} 
{\rm Let $Y$ be an arbitrary smooth manifold with boundary, $\nu\to Y$ a (twisted) spinor bundle and $\mu\to\partial Y$ a vector bundle. A  pseudo--differential operator $\mb{B}:C^{\infty}(\partial Y,\nu)\to C^{\infty}(\partial Y,\mu)$ of order $0$ is said to define an {\em elliptic boundary condition} (abbreviated e.b.c.) if and only if
\begin{enumerate}
	\item\label{extension} the extension $\mb{B}^{(s)}:H^s(\partial Y,\nu)\to H^s(\partial Y,\mu)$ to the chain of Sobolev spaces $H^s(\partial Y,\nu)$ and $H^s(\partial Y,\mu)$, $s\geq0$, has closed range.
	\item the restriction of the principal symbol $b=\sigma(\mb{B})_{|\mathrm{range}\big(\sigma(\mb{Q}_{\mb{D}})\big)}:\mathrm{range}\big(\sigma(\mb{Q}_{\mb{D}})\big)\to\mathrm{range}(b)$ is an isomorphism. 
\end{enumerate}
Furthermore, an e.b.c.\ is said to be {\em local}, if in addition $\mathrm{range}(p,\xi)=\nu_p$ holds for all $p\in\partial Y$ (in this case,~\ref{extension}. is automatically satisfied, cf.~{\rm\cite{bowo93}}, Rem.\ 18.2).} 
\end{definition}

The following theorem summarises the main properties of local e.b.c..

\begin{thm}\label{ebcprop}\hspace{-4pt}{\rm (\cite{bowo93}, Thm.\ 19.1, Thm.\ 20.12, Thm.\ 20.13 and Thm.\ 21.5)}
If $\mb{B}$ defines a local e.b.c., then 
\begin{enumerate}
	\item regularity holds, that is, $s\in H^s(Y,\nu)\cap\ker(\mb{D}\oplus\mb{B})$ implies $s\in C^{\infty}(Y,\nu)$.
	\item the operator $\mb{D}\oplus\mb{B}$ is Fredholm with index
$$
\ind(\mb{D}\oplus\mb{B})=\ind\big(\mb{B}\mb{Q}_{\mb{D}}:\mc{H}(\mb{D})\to C^{\infty}(\partial Y,\nu)\big).
$$
	\item $\ind(\mb{D}\oplus\mb{B})$ depends only on the homotopy type of the principal symbols involved. 
\end{enumerate}
Furthermore, if $Y$ is odd--dimensional, the orthogonal projector $\mb{P}^{\pm}$ onto the space of positive and negative half--spinors $\nu^{\pm}$ over the (even--dimensional) boundary always defines a local e.b.c.\ with $\ind(\mb{D}\oplus\mb{P}^{\pm})=0$. 
\end{thm}

\begin{rmk}
For even--dimensional manifolds the existence of local e.b.c.\ is topologically obstructed (e.g.~\cite{blbo85}, Section II.7.B).
\end{rmk}

The following result is a valuable tool for explicit index computations.

\begin{prp}\hspace{-4pt}{\rm (\cite{bowo93}, Thm.\ 21.2)}\label{indextool}
Let $\mb{D}:C^{\infty}(Y,\nu)\to C^{\infty}(Y,\nu)$ be a Dirac operator on some (twisted) spinor bundle $\nu$ over an odd--dimensional manifold $Y$ with boundary. Further, consider two orthogonal projectors onto subbundles $\nu_{1,2}$ of $\nu_{|\partial Y}$, $\mb{B}_1:C^{\infty}(\partial Y,\nu)\to C^{\infty}(\partial Y,\nu_1)$ and $\mb{B}_2:C^{\infty}(\partial Y,\nu)\to C^{\infty}(\partial Y,\nu_2)$, and suppose they define a local e.b.c.. Then
$$
\ind(\mb{D}\oplus\mb{B}_2)-\ind(\mb{D}\oplus\mb{B}_1)=\ind\big(\mb{B}_2\mb{Q}_{\mb{D}}\mb{B}^*_1:C^{\infty}(\partial Y,\nu_1)\to C^{\infty}(\partial Y,\nu_2)\big).
$$
\end{prp}

Coming back to equation~(\ref{defprob}), we shall write $\ind(X,Y)$ for $\ind(\mb{D}\oplus\mb{B})$.  Let $\overline{\partial}_{\nu_X}$ denote the Cauchy--Riemann operator associated with the natural complex structure (or equivalently, the natural orientation) of $\nu_X$ (cf.\ Lemma~\ref{boundary}). Finally, we are in a position to prove the central theorem of this paper.

\begin{thm}\label{index}
The operator $\mb{B}$ in~(\ref{defprob}) defines a local e.b.c.\ with index
$$
\ind(X,Y)=\ind(\overline{\partial}_{\nu_X}).
$$
\end{thm}
\begin{prf}
Let us fix some collar neighbourhood $\mc{C}\cong\partial Y\times[0,\epsilon)$ of $\partial Y$ for which we may assume the riemannian structure to be a product (possibly after homotopically deforming the metric). Further, we complete the inward pointing coordinate vector $u$ to a local orthonormal basis $\big(v(y,t),\,w(y,t)\big)$ of $T_y\partial Y\times\{t\}$ such that $u\times v=w$.  Near the boundary, we have the decomposition $\mb{D}=u\times(\nabla_u+\mb{R})$ with $\nabla_u=\partial_u$ and
\begin{equation}\label{R}
\mb{R}=w\times\nabla_v-v\times\nabla_w,
\end{equation}
as follows from $(a\times b)\times c=-a\times (b\times c)$ valid whenever $\{a,b,c\}$ is an orthogonal family, cf.~(\ref{chirule}).  Note that the bundles of positive and negative half--spinors $\nu^{\pm}$ inside $\nu^{\C}$ are just the eigenspaces of $J=u\times$.

\bigskip
   
Locally, we will work with the following basis of $\nu^{\C}$: Choose a nowhere vanishing local section $a\in C^{\infty}(\partial Y,\nu_X)$ so that
$$
\nu_X^{\C}\loc\langle\alpha=a-iJa,\overline{\alpha}=a+iJa\rangle,
$$
cf.\ Lemma~\ref{boundary}. Consider then the sections $\beta=-v\times\overline{\alpha}$ and $\overline{\beta}=-v\times\alpha$. Again, the lemma implies that
$$
\mu^{\C}_X\loc\langle\beta,\overline\beta\rangle.
$$
Furthermore,
$$
\nu^+\loc\langle\alpha,\beta\rangle,\quad\nu^-\loc\langle\overline{\alpha},\overline{\beta}\rangle.
$$
As an example, take $J\alpha=Ja+ia=i\alpha$ and $J\overline{\beta}=v\times(u\times\alpha)=-i\overline{\beta}$ etc.. For any subsequent matrix representation over $\nu^{\C}$, the ordered basis $\{\alpha, \beta, \overline{\alpha},\overline{\beta}\}$ shall be used.

\bigskip

Checking that $\mb{B}$ defines a local e.b.c.\ requires the principal symbol $q=\sigma(\mb{Q}_{\mb{D}})$ of the associated Calder\'on operator. By the Calder\'on--Seeley theorem (cf.\ Thm.\ 12.4 in~\cite{bowo93}), $q$ is  the projector onto the eigenspace of $\sigma(\mb{R})$ corresponding to eigenvalues with positive real part. With respect to our fixed local ordered basis of $\nu^{\C}$ around $p\in\partial Y$, $v$ and $w$ act as
$$
v\times=\left(\begin{array}{cc}\mbox{\hspace{-8pt}}\mb{0}&\mbox{\hspace{-8pt}}\begin{array}{cc}\scriptstyle\phantom{-}0&\scriptstyle\phantom{-}1\\\scriptstyle-1&\scriptstyle\phantom{-}0\end{array}\\\mbox{\hspace{-8pt}}\begin{array}{cc}\scriptstyle\phantom{-}0&\scriptstyle\phantom{-}1\\\scriptstyle-1&\scriptstyle\phantom{-}0\end{array}&\mb{0}\end{array}\right),\quad
w\times=\left(\begin{array}{cc}\mbox{\hspace{-8pt}}\mb{0}&\mbox{\hspace{-8pt}}\begin{array}{cc}\scriptstyle\phantom{-}0&\scriptstyle-i\\\scriptstyle\phantom{-}i&\scriptstyle\phantom{-}0\end{array}\\\mbox{\hspace{-8pt}}\begin{array}{cc}\scriptstyle\phantom{-}0&\scriptstyle\phantom{-}i\\\scriptstyle-i&\scriptstyle\phantom{-}0\end{array}&\mb{0}\end{array}\right).
$$
This follows from $v\times\alpha=-\overline{\beta}$, $w\times\alpha=-u\times(v\times\alpha)=-i\overline{\beta}$ etc.. For $(\eta_v,\eta_w)\in T^*_p\partial Y\backslash\{0\}$ of unit norm, we deduce from~(\ref{R}) (with $\eta=\eta_v+i\eta_w$) that
$$
\sigma(\mb{R})(p,\eta)=i(\eta_v\cdot w\times-\eta_w\cdot v\times)=\left(\begin{array}{cc}\mbox{\hspace{-8pt}}\mb{0}&\mbox{\hspace{-8pt}}\begin{array}{cc}\scriptstyle\phantom{-}0&\scriptstyle\phantom{-}\overline{\eta}\\\scriptstyle-\overline{\eta}&\scriptstyle\phantom{-}0\end{array}\\\mbox{\hspace{-8pt}}\begin{array}{cc}\scriptstyle\phantom{-}0&\scriptstyle-\eta\\\scriptstyle\phantom{-}\eta&\scriptstyle\phantom{-}0\end{array}&\mb{0}\end{array}\right)=\left(\begin{array}{cc}\mb{0}& r_-(p,\eta)\\r_+(x,\eta)&\mb{0}\end{array}\right).
$$ 
Now $r_+(p,\eta)^*=r_-(p,\eta)$ and $r_+(p,\eta)=r_-(p,\eta)^{-1}$, so that $\sigma(\mb{R})(p,\eta)^*=\sigma(\mb{R})(p,\eta)=\sigma(\mb{R})(p,\eta)^{-1}$. Consequently, the eigenvalues are $\pm1$, and the projector on the eigenspace associated with $1$ is given by
$$
q(p,\eta)=\frac{1}{2}\left(\begin{array}{cc}{\rm Id}_2&r_-(p,\eta)\\r_+(p,\eta)&{\rm Id}_2\end{array}\right).
$$
On the other hand, $\mb{B}$ is the orthogonal projector onto $\mu_X^\C$, so that its principal symbol is the matrix (taken with respect to the fixed basis of $\nu^{\C}$ and $\{\beta,\overline{\beta}\}$ of $\mu^{\C}_X$)
$$ 
\sigma(\mb{B})(p,\eta)=\left(\begin{array}{cccc} 0 & 1 & 0 & 0\\ 0 & 0 & 0 &1\end{array}\right),
$$
which is of full rank. Since
$$
\sigma(\mb{B})\circ q(p,\eta)=\frac{1}{2}\left(\begin{array}{cccc} 0 & 1 & -\overline{\eta} & 0\\ \eta & 0 & 0 &1\end{array}\right)
$$
is also of full rank, the boundary condition defined by $\mb{B}$ is local elliptic according to Definition~\ref{locbond}.

\bigskip

It remains to compute the index. As in Theorem~\ref{ebcprop}, let $\mb{P}^+$ denote the orthogonal projector onto $\nu^+$. We recall that according to this theorem, $\mb{P}^+$ defines a local e.b.c.\ with vanishing index. In virtue of Theorem~\ref{indextool} and the established local ellipticity of $\mb{B}$,
$$
\ind(X,Y)=\ind(\mb{D}\oplus\mb{B}) = \ind\big(\mb{B}\mb{Q}_{\mb{D}}\mb{P}^+:C^{\infty}(\partial Y,\nu^+) \to C^{\infty}(\partial Y,\mu_X^\C)\big).
$$
But the symbol of $\mb{B}\mb{Q}_{\mb{D}}\mb{P}^+$ is just
\begin{eqnarray*}
\sigma(\mb{B}\mb{Q}_{\mb{D}}\mb{P}^+)(p,\eta) & = & \sigma(\mb{B})\circ q\circ\sigma(\mb{P}^+)(p,\eta)\\
& = & \frac{1}{2}\left(\begin{array}{cccc} 0 & 1 & -\overline{\eta} & 0\\ \eta & 0 & 0 &1\end{array}\right)\left(\begin{array}{cc}1 & 0\\0 & 1\\0 & 0\\0 & 0\end{array}\right)\\
& = & \frac{1}{2}\left(\begin{array}{cc} 0 & 1\\\eta & 0\end{array}\right):\nu^+_p\to(\mu_X^{\C})_p,
\end{eqnarray*}
where the matrix is taken with respect to the basis $\{\alpha,\beta\}$ of $S^+$ and $\{\beta,\overline{\beta}\}$ of $\mu^{\C}_X$. In particular, the symbol sends $\beta$ to $\beta$ and therefore acts as the identity on $\mu^{1,0}_X=\nu^+\cap\mu^{\C}_X$. On the other hand, the induced map $\nu^{1,0}_X=\nu^+\cap\nu^{\C}_X\to \mu^{0,1}_X=\nu^-\cap\mu^{\C}_X$ is up to $-i$ the symbol of the Cauchy--Riemann operator $\overline{\partial}_{\nu_X}$ on $\nu_X^{1,0}$ after the identification $\mu^{0,1}_X\cong\nu^{1,0}_X\otimes\overline{K}_{\partial Y}$ (cf.\ Lemma~\ref{boundary}). Indeed, on a trivialisation of $\nu^{\C}_X$ it acts as $\overline{\partial}_{\nu_X}=(\partial_1+i\partial_2)/2$, where $(x_1,x_2)$ are coordinates such that $\partial_1(p)=v(p)$ and $\partial_2(p)=w(p)$. Hence $\sigma(\overline{\partial}_{\nu_X})(p,\eta)=i(\eta_v+i\eta_w)/2$ which is what we wanted.
\end{prf}

As a consequence of Riemann--Roch, we obtain the

\begin{cor}
Let $\Sigma_{g_j}$ be the connected components of $\partial Y$ of genus $g_j$, and $c_1(\nu_{X|\Sigma_{g_j}})$ be the first Chern classes of $\nu_{X|\Sigma_{g_j}}$. Then
$$
\ind(X,Y)=\sum_j\int_{\Sigma_{g_j}}c_1(\nu_{X|\Sigma_{g_j}})+1-g_j.
$$
\end{cor}
%
%
%
%
%
%
\section{$\psi$--positive boundary conditions}
\label{gendefthm}
For the deformation problem~(\ref{defprob}) the bundle $\nu_X$ is the only non--intrinsic piece of datum attached to $Y$, and its properties were derived using the coassociativity of $X$. In view of the previous theorem, we may relax the coassociativity condition as follows:

\begin{definition}
{\rm Let $(M,\varphi)$ be a $G_2$--manifold and $\psi=\star\varphi$. An orientable $4$--submanifold of $M$ is said to be $\psi$--{\em positive at} $p\in X$ if and only if its tangent space at $p$ is $\psi$--positive for some orientation, i.e.\ $\psi_{p|T_pX}>0$. Further, $X$ is called $\psi$--{\em positive} if for a suitable orientation $X$ is $\psi$--positive at all points $p\in X$.}
\end{definition}

Obvious examples of $\psi$--positive submanifolds are coassociatives. An alternative characterisation of $\psi$--positivity is the notion of $\varphi$--freeness going back to recent work of Harvey and Lawson~\cite{hala06}.

\begin{lem}\label{psivsphi}
A suitably oriented $4$--plane $F\in G_4(\R^7)$ is $\psi$--positive at $p\in X$ if and only if $F$ is $\varphi${\em--free}, i.e.\ $F$ contains no associative $3$--plane. 
\end{lem}
\begin{prf}
If $F\in G_4(T_xM)$ is not $\varphi$--free, then $F$ contains an associative subplane $E$. Hence $\varphi_{|E^{\perp}}\equiv 0$ and in particular $\varphi_{|F^{\perp}}=\psi_{|F}\equiv 0$, which proves necessity. Conversely, assume that $\psi_{|F}\equiv0$. We write $F^{\perp}=u\wedge v\wedge w$ and choose a vector $a$ in the orthogonal complement to the linear span of $u,\,v,\,w,\,u\times v,\,u\times w,\,v\times w$, which at most is six dimensional. As a result, the $4$--plane $E^{\perp}$ spanned by $F^{\perp}$ and $a$ is coassociative, for $\varphi_{|E^{\perp}}\equiv0$. Hence $E\subset F$ is associative.
\end{prf}

Put differently, for any $a,\,b\in T_pX$, we have $a\times b\notin T_pX$ if $X$ is $\psi$--positive at $p$. From this point of view the class of $\psi$--positives inside $G_2$--manifolds naturally matches the class of totally real submanifolds inside K\"ahler manifolds (which, in particular, comprises lagrangians). 

\bigskip

Next we wish to investigate the Zariski tangent space of $\mf{M}_{X,Y}$ under the assumption that $X$ is merely $\psi$--positive. To that end, let $N_X\to\partial Y$ be the orthogonal complement of $T\partial Y$ in $TX_{|\partial Y}$, and define $\nu_X$ to be the image of $N_X$ under the orthogonal projection $\mb{P}:TM_{|Y}\to\nu_{|Y}$. As before, $\mu_X\to\partial Y$ denotes the orthogonal complement of $\nu_X$ in $\nu_{|\partial Y}$. The geometry on the boundary is specified by the following

\begin{lem}\label{boundaryvarphifree}
If $X$ is $\psi$--positive, then
\begin{enumerate}
	\item\label{isoNnu} the restriction of  $\mb{P}$ to $N_X$ defines an isomorphism onto $\nu_X=\mb{P}(N_X)$. 
	\item for any non--zero $b\in\mu_X$, $Jb=u\times b\notin\nu_X$, that is for the orthogonal projection $\mb{p}_{\mu_X}$ on $\mu_X$ we have $\mb{p}_{\mu_X}(Jb)\not=0$. 
\end{enumerate}
\end{lem}
\begin{prf} 
The kernel of $\mb{P}_{|TM_{|\partial Y}}$ is $TY_{|\partial Y}$. Since $TX$ is $\psi$--positive, $N_X\cap\ker\pi=\{0\}$ according to the previous proposition, whence the first assertion. Next, suppose there is a $b_0\in\mu_X$ of unit norm such that $Jb_0\in\nu_X$. Let $b_1\in\mu_X$ be a vector orthogonal to $b_0$. Then $Jb_1$ lies in $\nu_X$, for it is orthogonal to $b_1$ and $g(b_0,Jb_1)=-g(Jb_0,b_1)=0$. By~\ref{isoNnu}.\ there exist uniquely determined $n_{0,1}\in N_X$ with $\mb{P}(n_{0,1})=Jb_{0,1}$. Since $N_X\perp T\partial Y$, we have in fact $n_{0,1}=Jb_{0,1}+\lambda_{0,1}u$ for $\lambda_{0,1}\in\R$ (where $u$ denotes again the inward pointing normal vector on the boundary). It is straightforward to check that the orthogonal complement of $T_{\pi(b_0)}X$ in $T_{\pi(b_0)}M$ is spanned by $A=b_0$, $B=b_1$ and $C=u-\lambda_0Jb_0-\lambda_1Jb_1$. Moreover, $v=A\times B\in TY_{|\partial Y}$ belongs in fact to $T\partial Y$, for $g(u,v)=\varphi(b_0,b_1,u)=-g(Jb_1,b_0)=0$. Consequently, $\varphi(A,B,C)= g(A\times B,C)=0$ which contradicts the $\psi$--positiveness of $TX$. 
\end{prf}

Consequently, $N=N_X\oplus\mu_X$ and $\nu_X\oplus\mu_X$ are isomorphic subbundles of $TM_{|\partial Y}$. Further, we can extend $N\to\partial Y$ to a new subbundle $N\subset TM_{|Y}$ transversal to $TY$, and adapt Corollary~\ref{defequs} to this more general situation: If $N_X$ is orthogonal to $TY_{|\partial Y}$ (this does not imply the fibres of $N$ to be coassociative), we can take $N=\nu$ and identify $T^{\mathrm{Zar}}\mc{M}_{X,Y}$ with the space of solutions of~(\ref{defprob}). Otherwise first extend $N_X$ trivially on some collar neighbourhood. Since being isomorphic is a stable property, we can homotope $N_X$ over $\mc{C}$ into a new bundle (still written $N_X$) with fibres in $\nu_X$ sufficiently far away from the boundary. Then extend $N=N_X\oplus\mu_X$ by $\nu$. A section $s\in C^{\infty}(Y,N)$ induces a family of submanifolds $Y_t$ with boundary in $X$ if $s_{|\partial Y}$ lies in the kernel of $\mb{B}: C^{\infty}(\partial Y,N)\to C^{\infty}(\partial Y,\mu_X)$, the orthogonal projection onto $\mu_X$ in $N_{|\partial Y}$. Further, for $Y_t$ to be associative we need the normal component $\mb{P}(s)$ of $s$ to be harmonic. Thus $T^{\mathrm{Zar}}\mc{M}_{X,Y}$ can be identified with the solution space of
\begin{equation}\label{extdefprob}
\mb{D}\mb{P}(s)=0,\quad\mb{B}(s_{|\partial Y})=0,\quad s\in C^{\infty}(Y,N).
\end{equation}

\begin{prp}
Let $Y\subset M$ be an associative with boundary inside a $\psi$--positive submanifold $X$. Then 
$$
\ind(X,Y)=\ind(\overline{\partial}_{\nu_X}).
$$ 
\end{prp}

\begin{rmk}
Here, the Cauchy--Riemann operator on $\nu_X$ is defined by the complex structure coming from the induced orientation on $\nu_X\cong N_X$.
\end{rmk}

\begin{prf} 
For simplicity we assume $\mb{P}=\mathrm{Id}$, i.e.\ $N_X$ is orthogonal to $TY_{|\partial Y}$. We choose a nowhere vanishing local section $a$ of $\nu_X$ and extend $b=-v\times a$ to a local trivialisation $\{b,\widetilde{b}\}$ of $\mu_X$. By the previous lemma we may take $\widetilde{b}$ to be the orthogonal projection of $Jb$ to $\mu_X$. Let $(0,1,s,t)$ be the coordinates of $\widetilde{b}$ with respect to the local basis $\{b,Jb,v\times b,w\times b\}$ of $\nu$, where $v$ is a nowhere vanishing local section of $T\partial Y$ and $w=Jv=u\times v$. The latter basis gives rise to the basis $\{\alpha,\beta,\overline{\alpha},\overline{\beta}\}$ of $\nu^{\C}$ as given in the proof of Theorem~\ref{index} with respect to which the matrix of $\sigma(\mb{B}^{\C})$ can be written as
$$ 
\sigma(\mb{B}^\C)(x,\eta)=\left(\begin{array}{cccc} 0 & 1 & 0 & 1\\ z & -i & \overline{z} &i\end{array}\right),
$$
where $z=s+it$. For $\mb{B}^\C\mb{Q}\mb{P}^+$, we therefore find
\begin{equation}\label{defps}
\sigma(\mb{B}^\C\mb{Q}\mb{P}^+)(x,\eta)=\frac{1}{2}\left(\begin{array}{cc} \eta & 1\\z+i\eta & -i-z\eta\end{array}\right):S_x^+\to(\mu_X^{\C})_x.
\end{equation}
The determinant of this matrix is $(-2i\eta-z-\overline{z}\eta^2)/4$, and multiplication with $\overline{\eta}$ shows this to vanish only if
$\Re(\overline{T}\eta)=-i$. Hence the system~(\ref{extdefprob}) is still elliptic. Furthermore, $F(Jb,t)=t\mb{p}_{\nu_X}(Jb)\oplus\mb{p}_{\mu_X}(Jb)$ is a global homotopy deforming $J\mu_X$ into $\mu_X$ and which in particular deforms $Jb$ into $\widetilde{b}$. Consequently, the symbol~(\ref{defps}) is homotoped into the symbol of $\overline{\partial}_{\nu_X}$. By homotopy invariance we recover the same index as before.
\end{prf}
%
%
%
%
%
%
%
%
\section{Examples}
\label{indexa}
In this section we construct explicit examples of pairs $(X,Y)$ with $\ind(X,Y)\not=0$.  Throughout this section we denote by $\Sigma_g$ a surface diffeomorphic to a compact oriented Riemann surface of genus $g$.

\bigskip

\textbf{Local coassociative submanifolds.}
The first example will be local in nature, that is, the boundary of the compact associative will be contained in a local coassociative submanifold. Existence will be established by using Cartan--K\"ahler theory which requires all geometric objects involved (manifolds, boundaries, maps etc.) to be real analytic. Note that a torsion--free $G_2$--manifold $(M,g)$ is Ricci--flat so that the underlying riemannian metric is real analytic in harmonic coordinates~\cite{dk81}. Since $\varphi$ is harmonic \wrt $g$ (cf.\ Section~\ref{g2manif}), it will be real analytic in these coordinates, too. Consider then an associative $Y$  with real analytic boundary $\partial Y$, and a nowhere vanishing real analytic section $a\in C^{\infty}(\partial Y,\nu)$. The geodesic flow $\gamma_a:\partial Y\times(-\epsilon,\epsilon)\to M$ induced by $a$ is also real analytic in harmonic coordinates and therefore generates an analytical submanifold $N$ of dimension $3$. Further, $\varphi(v,w,a)=0$ for $v,w\in T\partial Y$. We conclude that the pull--back of $\varphi$ to $N$ vanishes identically, for $\nabla\varphi=0$. A Cartan--K\"ahler type argument invoked by Harvey and Lawson~\cite{hala82} (see also~\cite{br00}) shows that $N$ is contained in a local coassociative $X$. Further, $\nu_X$ is generated by $a$ and $u\times a$, where $u$ denotes again the inward pointing normal vector field of $\partial Y$. Thus $c_1(\nu_X)=0$, whence $\ind(X,Y)=1-g$. For example, taking $M=\Im\Hh\oplus\Hh$ and $\partial Y=\Sigma_g\subset\Im\Hh$, where $\Sigma_g$ is a $g$--handle body in $\Im\Hh\cong\R^3$ and $Y$ the relatively compact interior bounded by $\Sigma_g$, yields examples of arbitrary negative index. An actual deformation is given by moving $Y$ along the straight line determined by $a$. In fact, taking for $\Sigma_0$ the standard $2$--sphere bounding the unit ball, one obtains a smooth moduli space which is thus of actual dimension $1$ (cf.~\cite{ga09}).

\bigskip

\textbf{Local associative submanifolds.} 
In the same vein, consider a real analytic surface $\Sigma_g$ inside some real analytic coassociative $X$ in $M$. Again, Cartan--K\"ahler theory yields the existence of a local associative $Y$ containing $\Sigma_g$~\cite{hala82},~\cite{rosa07}. Using a collar neighbourhood of $\Sigma_g$ inside $Y$ we can construct an associative of the form $\Sigma_g\times [0,1]$ which we keep on denoting by $Y$ for simplicity. Further, we can translate $X$ into a submanifold $X'$ containing $\Sigma_g\times\{1\}$ by a suitable diffeomorphism $C^1$--close enough to the identity. Of course, there is no reason for this diffeomorphism to preserve the $G_2$--structure, so that $X'$ will not be coassociative in general. However, as $\psi$--positivity is a pointwise open condition, $X'$ will be still $\psi$--positive for suitable choices, and the generalised deformation theorem applies. To compute the index we note that $N_{X'}\to\Sigma_g\times\{1\}$ is obtained by translating $\nu_X\to\Sigma_g\times\{0\}$, but the orientation flips as $u$ which points inward at $\Sigma\times\{0\}$ points outward at $\Sigma_g\times\{1\}$. Hence $c_1(\nu_X)=-c_1(\nu_{X'})$. Since the orientation also flips on $\Sigma_g\times\{1\}$ we finally get $\ind(X\cup X',Y)=2\ind(\overline{\partial}_{\Sigma_g})$.

\bigskip

{\bf The Bryant--Salamon construction.}
In~\cite{brsa87} Bryant and Salamon constructed holonomy $G_2$--metrics on (an open set of) the total space of the spinor bundle $S\to M^3$, where $M$ is a three--dimensional space form. In particular, when $M$ is taken to be the round sphere $\mc{S}^3$, there exists a {\em complete} holonomy $G_2$--metric on the total space $S\cong\mc{S}^3\times\Hh$ such that the fibres are orthogonal to the horizontal distribution of the canonical spin connection induced by $\nabla^g$. Furthermore, the zero section defines an associative. Any 3--ball $Y$ inside $\mc{S}^3$ is therefore associative with boundary $\partial Y=\mc{S}^2=\Sigma_0$. Let $u$ be again the inward pointing normal vector field near the boundary. Assume that $S_{|\Sigma_0}$ has a subbundle which is a $J=u\times$--complex line bundle over $\Sigma_0$ of degree $n$. If $X_n$ denotes the total space, then $X_n$ is coassociative at $\partial Y$ (cf.\ Lemma~\ref{boundary}), hence $\psi$--positive near $\partial Y$. Therefore, the index formula applies and yields $\ind(X_n,Y)=n+1$. For $n=0$, such a line bundle can be constructed by taking a nowhere vanishing section $a\in\Gamma(S_{|\Sigma_0})$. This gives the trivial complex line bundle spanned by $a$ and $u\times a$, whence $\ind(X_0,Y)=1$.

\bigskip

{\bf The Calabi--Yau extension.} 
Let $(K,\omega,\Omega)$ be a Calabi--Yau $3$--fold with K\"ahler $2$--form $\omega$ and holomorphic volume form $\Omega$. Then $\omega^l/l!$ and $\cos(\theta)\Re\Omega+\sin(\theta)\Im\Omega$ define calibrations on $K$ which calibrate complex submanifolds (of {\em real} dimension $2l$) and {\em special lagrangians of phase $e^{i\theta}$} respectively. On $M=K\times \mc{S}^1$ we can define a torsion--free $G_2$--structure by
$$
\varphi=\Re\Omega+\omega\wedge dt.
$$ 
Further, if $C\subset K$ is a complex curve or $L\subset K$ a special lagrangian of phase $1$, then $C\times \mc{S}^1$ and $L\times\{pt\}$ are associative. If $S\subset K$ is a complex surface or $L\subset K$ is a special lagrangian of phase $i$, then $S\times\{pt\}$ and $L\times \mc{S}^1$ are coassociative. Therefore, let $L$ be a special lagrangian of phase $1$ with boundary $\partial L=\Sigma_g$ inside a complex surface $S$. The normal bundle of $L$ is just $J(TL)$, for $L$ is lagrangian. On the other hand, $S$ is complex, so that $J(T\partial L)\subset TS$. Lifting $S$ and $L$ to the coassociative $X= S \times \{pt\}$ and the associative $Y=L \times \{pt\}$ inside $M$, we find that the underlying real rank $2$ bundle of $\nu_X$ coincides with $J(T\partial L)\oplus 0$. Now the normal bundle $\nu\to Y$ is generated by $J(TL)$ and $\partial_t$. Consequently, $\partial_t$ is a section of $\mu_X$, the orthogonal complement of $\nu_X$ inside $\nu$. Lemma~\ref{boundary} then implies $\nu_X\cong\overline{T\partial L}=K_{\partial L}$ as complex line bundles so that $c_1(\nu_X)=c_1(K_{\partial L})=2(1-g)$. We finally obtain $\ind(X,Y)=3(1-g)$ from the index formula. These considerations apply for instance to real, $3$--dimensional submanifolds contained in the real part of a smooth quintic inside $\CP^4$, and whose boundaries are real analytic. The real part condition guarantees that such submanifolds are special lagrangian, while the boundary condition ensures the boundary to be contained in a complex surface $S$.

\bigskip

{\bf Flat compact examples.}
On $\R^7$ we have the trivial $G_2$--structure $\varphi_0=dx^{123}+dx^1\wedge(dx^{45}+dx^{67})+dx^2\wedge(dx^{46}-dx^{57})+dx^3\wedge(-dx^{47}-dx^{56})$. Since $\varphi_0$ is translation invariant, it descends to a torsion--free $G_2$--structure on the torus $T^7=\R^7/\Z^7$ still denoted $\varphi_0$. Let $(x_1,\ldots,x_7)$, $x_i\in\R/\Z$ be coordinates on $T^7$. The isometric involutions
$$
\begin{array}{lcl}
\sigma_0(x_1,\ldots,x_7) & = & (x_1,x_2,x_3,-x_4,-x_5,-x_6,-x_7)\\
\tau_0(x_1,\ldots,x_7) & = & (-x_1,x_2,x_3,x_4,x_5,-x_6,-x_7)\\
\end{array}
$$
satisfy $\sigma_0^*\varphi_0=\varphi_0$ and $\tau_0^*\varphi_0=-\varphi_0$. In the sequel, let $\alpha_i$, $i=1,\ldots,7$ denote an element in $\{0,1/2\}$. The fixed point locus of $\sigma_0$ consists of the sixteen $3$--tori $\sigma_0 T^3_{...a_4a_5a_6a_7}=\{(x_1,x_2,x_3,a_4,a_5,a_6,a_7)\,|\,x_i\in\R/\Z\}$ which are associative as can be seen by direct inspection. Further, the fixed point locus of $\tau_0$ consists of the eight coassociative $4$--tori $\tau_0 T^3_{a_1....a_6a_7}=\{(a_1,x_2,x_3,x_4,x_5,a_6,a_7)\,|\,x_i\in\R/\Z\}$. Then $Y=\{(t,x_2,x_3,0,0,0,0)\,|\,0\leq t\leq1/2,\,x_i\in\R/\Z\}\cong[0,1/2]\times\Sigma_1$ is an associative with boundary in the coassociative $X=\tau_0 T^4_{0....00}\cup\tau_0 T^4_{\frac{1}{2}....00}$. Since $\nu_X$ is trivial, we find $\ind(X,Y)=0$. Of course, there exists actual deformations of $Y$ along $X$ induced by the flow of $\partial_{x_4}$ and $\partial_{x_5}$.

\bigskip

{\bf Joyce manifolds.}
For the remaining examples we first invoke the following result which generalises the previous example.

\begin{prp}\hspace{-4pt}{\rm (\cite{jo00}, Sec.\ 10.8)}
Let $(M,\varphi)$ be a torsion--free $G_2$--manifold.

(i) If $\sigma:M\to M$ is a nontrivial  isometric involution with $\sigma^*\varphi=\varphi$, then the fixed point locus of  $\sigma$ defines an associative submanifold. 

(ii) If $\tau:M\to M$ is a nontrivial  isometric involution with $\tau^*\varphi=-\varphi$, then each connected component of the fixed point locus of $\tau$ is either a coassociative 4--fold or a single point.
 \end{prp}

A method for the construction of compact holonomy $G_2$--manifolds is also due to Joyce. We briefly outline his method to fix notation (cf.~\cite{jo00} for details and~\cite{jo04} for a quick introduction). One considers quotients $T^7/\Gamma$ where $\Gamma$ is a discrete group acting as a subgroup of $\varphi_0$--preserving isometries of $(T^7,g_0)$. In particular, $\varphi_0$ descends to a $G_2$--form $\varphi$ on $T^7/\Gamma$. For instance choose the group $\Gamma\cong(\Z/2\Z)^3$ given in 12.2 of~\cite{jo00} which is generated by
$$
\begin{array}{lcl}
\alpha(x_1,\ldots,x_7) & = & (x_1,x_2,x_3,-x_4,-x_5,-x_6,-x_7)\\
\beta(x_1,\ldots,x_7) & = & (x_1,-x_2,-x_3,x_4,x_5,\frac{1}{2}-x_6,-x_7)\\
\gamma(x_1,\ldots,x_7) & = & (-x_1,x_2,-x_3,x_4,\frac{1}{2}-x_5,x_6,\frac{1}{2}-x_7).
\end{array}
$$
There are no fixed points of $\delta\in\Gamma\backslash\{\mathrm{Id}\}$ unless $\delta=\alpha,\,\beta$ or $\gamma$. For instance, $\alpha$ fixes the $3$--tori $\alpha T^3_{...a_4a_5a_6a_7}=\{(x_1,x_2,x_3,a_4,a_5,a_6,a_7)\,|\,x_i\in\R/\Z\}$. The subgroup $\langle\beta,\gamma\rangle$ acts freely on these sixteen tori. For example,  $\beta$ maps $\alpha T^3_{...0000}$ to $\alpha T^3_{...00\frac{1}{2}0}$ etc.. Hence, if we denote the image of these tori in $T^7/\Gamma$ by $\alpha T^3_{[...a_4a_5a_6a_7]}$, then $\alpha T^3_{[...0000]}=\alpha T^3_{[...00\frac{1}{2}0]}$ so that $\alpha$ contributes four distinct copies of $T^3$ to the singular locus of $T^7/\Gamma$. A similar analysis applies to $\beta$ and $\gamma$. In total, the singular locus consists of twelve distinct copies of $T^3$ which do not intersect. Near any singular $3$--torus we can fix a neighbourhood homeomorphic to $\R^3\times(\C^2/\{\pm1\})$ for instance by mapping $[x_1,x_2,x_3,a_4+y_4,a_5+y_5,a_6+y_6,a_7+y_7]$ near $\alpha T^3_{[...a_4a_5a_6a_7]}$ to $(x_1,x_4,x_5,\{z_1=y_2+iy_3,z_2=y_6+iy_7\})$ (with $\{z_1,z_2\}=\{-z_1,-z_2\}$ the coset of $(z_1,z_2)$ inside $\C^2/\{\pm1\}$). Blowing up the singular point in $\C^2/\{\pm1\}$ yields a resolution $R$ on which there exists a family $h_t$ of asymptotically flat Calabi--Yau metrics. Together with the flat $3$--torus $(T^3,g_0)$, the product metric $g_0\oplus h_t$ on $T^3\times R$ gives a torsion--free $G_2$--structure which can be glued into $T^7/\Gamma$. Proceeding this way with the remaining singular tori yields a smooth compact resolution $\pi:M\to T^7/\Gamma$ on which we have a family of $G_2$--structures $\widetilde{\varphi}_t$ with ``small'' torsion. For $t$ sufficiently small, Joyce's deformation theorem ensures the existence of a torsion--free $G_2$--structure $\widehat{\varphi}_t$ in a $C^0$--vicinity of $\widetilde{\varphi}_t$ controlled by $t$ which induces a holonomy $G_2$--metric for $M$ is simply--connected. 

\bigskip

As a warm--up we consider the example induced by the isometric involutions
$$
\begin{array}{lcl}
\sigma_0(x_1,\ldots,x_7) & = & (x_1,\frac{1}{2}-x_2,\frac{1}{2}-x_3,x_4,x_5,-x_6,\frac{1}{2}-x_7)\\
\tau_0(x_1,\ldots,x_7) & = & (x_1,x_2,\frac{1}{2}-x_3,\frac{1}{2}-x_4,x_5,x_6,\frac{1}{2}-x_7).
\end{array}
$$
of $(T^7,\varphi_0)$. Clearly, $\sigma^*_0\varphi_0=\varphi_0$ and $\tau^*_0\varphi_0=-\varphi_0$. Further, $\sigma_0$ and $\tau_0$ commute with $\Gamma$ and thus descend to isometric involutions $\sigma$ and $\tau$ of $T^7/\Gamma$. The fixed point set of $\sigma_0$ is given by the sixteen $3$--tori $\sigma_0T^3_{.b_2b_3..a_6b_7}=\{(x_1,b_2,b_3,x_4,x_5,a_6,b_7)\,|\,x_i\in \R/\Z\}$ where from now on, $b_j$ denotes an element in $\{1/4,3/4\}$. Furthermore, $\sigma_0\circ\delta$ has no fixed points for $\delta\in\Gamma$ unless $\delta=\mathrm{Id}$. As $\Gamma$ acts transitively on the set of the tori $\sigma_0T^3_{.b_2b_3..a_6b_7}$ the resulting fixed point set of $\sigma$ consists of two tori, namely $\sigma T^3_{[.\frac{1}{4}\frac{1}{4}..0\frac{1}{4}]}$ and $\sigma T^3_{[.\frac{1}{4}\frac{1}{4}..\frac{1}{2}\frac{1}{4}]}$ which do not hit the singular locus of $T^7/\Gamma$. Similarly, $\tau_0\circ\delta$ has no fixed points unless $\delta=\mathrm{Id}$ (giving eight $4$--tori $\tau_0T^4_{..b_3b_4..b_7}$), or $\delta=\beta\gamma$ (in which case the fixed point locus consists of the $128$ fixed points $(a_1,a_2,b_3,b_4,b_5,b_6,b_7)$). The fixed point locus of $\tau$ is thus $\tau T^4_{[..\frac{1}{4}\frac{1}{4}..\frac{1}{4}]}$ plus eight isolated points. Again, this does not hit the singular locus. Now as in 12.6 of~\cite{jo00} one can resolve $T^7/\Gamma$ in a $\sigma,\tau$--equivariant way such that $\sigma$ and $\tau$ lift to isometric involutions preserving and reversing the torsion--free $G_2$--form $\widehat{\varphi}_t$ respectively. Since the fixed point loci of $\sigma$ and $\tau$ are not modified by the resolution process, we obtain an associative $Y$ coming from $\{[x_1,\frac{1}{4},\frac{1}{4},t,x_5,0,\frac{1}{4}]\,|\,x_{1,5}\in\R/\Z,\,\frac{1}{4}\leq t\leq\frac{3}{4}\}\cong[1/4,3/4]\times\Sigma_1$ which intersects the coassociative $X$ coming from $T^4_{[..\frac{1}{4}\frac{1}{4}..\frac{1}{4}]}$ in the 2--tori coming from $T^2_{[.\frac{1}{4}\frac{1}{4}\frac{1}{4}.0.]}$ and $T^2_{[.\frac{1}{4}\frac{1}{4}\frac{3}{4}.0.]}$. In particular, $\nu_X$ is trivial so we find $\ind(X,Y)=2(1-g)=0$.

\bigskip

For the second example we start with the isometric involutions
$$
\begin{array}{lcl}
\sigma_0(x_1,\ldots,x_7) & = & (x_1,\frac{1}{2}-x_2,\frac{1}{2}-x_3,x_4,x_5,-x_6,-x_7)\\
\tau_0(x_1,\ldots,x_7) & = & (\frac{1}{2}-x_1,x_2,x_3,x_4,x_5,-x_6,-x_7)
\end{array}
$$
which also satisfy $\sigma_0^*\varphi_0=\varphi_0$ and $\tau_0^*\varphi_0=-\varphi_0$ respectively. Now $\sigma_0\circ\delta$ has no fixed points unless $\delta=\mathrm{Id}$ or $\alpha$. The fixed point locus of $\sigma$ is therefore the image under the quotient map of the sixteen tori $\sigma_0T^3_{.b_2b_3..a_6a_7}$ fixed by $\sigma_0$ and the sixteen tori $\sigma_0\alpha T^3_{.b_2b_3a_4a_5..}$ fixed by $\sigma_0\alpha$. Similarly, $\tau_0\circ\delta$ has no fixed points unless $\delta=\mathrm{Id}$ (in which case the fixed point locus consists of the eight tori $\tau_0T^4_{b_1....a_6a_7}$), $\delta=\alpha$ (giving eight tori $\tau_0\delta T^4_{b_1..a_4a_5..}$) or $\delta=\alpha\circ\beta$ \big(giving 128 points $(b_1,a_2,a_3,a_4,a_5,b_6,a_7)$\big). Consequently, the fixed point loci of both $\sigma$ and $\tau$ hit the singular $\alpha$--tori of $T^7/\Gamma$. To resolve $T^7/\Gamma$ in a $\sigma,\,\tau$--equivariant way and to compute their fixed point locus, we therefore proceed in two steps (cf.\ also 12.6.3 in~\cite{jo00}). Firstly, we consider the induced involutions on $T^7/\langle\alpha\rangle\cong T^3\times(T^4/\langle\pm1\rangle)$ also denoted by $\sigma$ and $\tau$, and resolve in a $\sigma,\,\tau$--equivariant way. The result is diffeomorphic to $T^3\times K3$, where $K3$ denotes the $K3$ surface. To resolve $T^7/\Gamma$ completely, we lift the actions of $\beta$ and $\gamma$ to $T^3\times K3$ and resolve the orbifold $(T^3\times K3)/\langle\beta,\gamma\rangle$ in a $\sigma,\,\tau$--equivariant way. First the fixed point locus of $\sigma$ inside $T^7/\langle\alpha\rangle$ consists of the sets $T^1_{.b_2b_3}\times S\cong\{(x_1,b_2,b_3)\,|\,x_1\in\R/\Z\}\times S$. Here, $S$ is the singular connected surface given by the union of $T^2_{..a_6a_7}/\langle\pm1\rangle\cong\{(x_4,x_5,a_6,a_7)\,|\,x_i\in\R/\Z\}/\langle\alpha\rangle$ with $T^2_{a_4a_5..}/\langle\pm1\rangle\cong\{(a_4,a_5,x_6,x_7)\,|\,x_i\in\R/\Z\}/\langle\alpha\rangle$ which intersect at the sixteen singular points $(a_4,\ldots,a_7)$. Similarly, the fixed point set of $\tau$ inside $T^7/\Gamma$ is $T^2_{b_1..}\times S=\{(b_1,x_2,x_3)\}\times S$. The lift of $\sigma$ and $\tau$ to the first resolution yields $T^1_{.b_2b_3}\times\Sigma$ and $T^2_{b_1..}\times\Sigma$ for some compact oriented Riemann surface $\Sigma$ as the corresponding fixed point set. Working out the zeroth and first Betti number of $\Sigma$ as in 12.5.2 of~\cite{jo00} shows that $b_0(\Sigma)=8$ and $b_1(\Sigma)=0$, that is, $\Sigma$ consists of eight copies of the $2$--sphere $\Sigma_0=\mc{S}^2$. In essence, every $T^2_{a_4a_5..}/\langle\pm1\rangle$ or $T^2_{..a_6a_7}/\langle\pm1\rangle$ contributes one $2$--sphere on the resolution which we label accordingly by $\mc{S}^2_{a_4a_5..}$ and $\mc{S}^2_{..a_6a_7}$. Now the action of $\beta$ and $\gamma$ on the fixed 32 copies of $T^1\times \mc{S}^2$ results in eight distinct copies inside $(T^3\times K3)/\langle\beta,\gamma\rangle$. For $T^2_{b_1..}\times\Sigma$ we find four distinct copies of $(T^2\times\mc{S}^2)/\langle\beta\rangle$ and two distinct copies of $T^2\times\mc{S}^2$, for $\gamma$ acts freely on the set of $2$--tori $T^2_{b_1..}$, while $\beta$ acts trivially and maps for instance $\mc{S}^2_{a_4a_5..}$ into itself, but $\mc{S}^2_{..a_6a_7}$ to $\mc{S}^2_{..\bar{a}_6a_7}$ etc.. Consequently, the boundary of the associative $Y$ coming from $\{[t,b_2,b_3]\,|\,t\in[1/4,3/4]\}\times\mc{S}^2_{[..a_6a_7]}$ inside $(T^3\times K3)/\langle\beta,\gamma\rangle$ consists of two $2$--spheres inside the coassociative $X$ coming from the union $X_{\frac{1}{4}}\cup X_{\frac{3}{4}}=T^2_{[\frac{1}{4}..]}\times\mc{S}^2_{[..a_6a_7]}\cup T^2_{[\frac{3}{4}..]}\times\mc{S}^2_{[..a_6a_7]}$. Inspection of the construction and the resolution process shows that $\nu_{X_{\frac{1}{4}}}$ and $\nu_{X_{\frac{3}{4}}}$ are trivial whence $\ind(X,Y)=2$.
%
%
%
%
%
\section{A $G_2$ analogue of the Maslov index}
Finally, we wish to introduce a $G_2$ analogue of the Maslov index. This index plays an important r\^ole in symplectic geometry, in particular for Floer homology, where it arises as a sort of relative Morse index~\cite{dosa94}, \cite{mcsa94},~\cite{saze92}.

\bigskip

Let us briefly recall its construction. We consider an almost complex manifold $(M^{2m},J)$ with an embedded (not necessarily holomorphic) $2$--disk $D$ with boundary in a totally real oriented submanifold $X^m$ (i.e for all $p\in X$, $T_pX$ does not contain any $J$--complex line). Since $D$ is contractible,  we can trivialise $TM_{|D}$ and regard the subbundle $TX_{|\partial D}$ as a closed curve in the set of totally real oriented $m$--planes in $\C^m$. On the other hand, this set is parametrised by\footnote{We denote by $GL_m^+(\R)$ the set of invertible $m\times m$--matrices with positive determinant.} $GL_m(\C)/GL^+_m(\R)$ and is homotopy equivalent to $U(m)/SO(m)$ -- the set of oriented lagrangian $m$--planes inside $\C^m$. By the exact homotopy sequence for fibrations $\pi_1\big(U(m)/SO(m)\big)\cong\Z$. The {\em Maslov index} $\mu(X,D)$ of $D$ is the integer corresponding to the homotopy class induced by $TX_{|\partial D}$.

\bigskip

The natural counterpart in the $G_2$--setting should be the following. Let $Y$ be an embedded (not necessarily associative) $3$--disk inside a topological $G_2$--manifold $M$ such that $\partial Y\cong\mc{S}^2$ lies in some $\psi$--positive submanifold $X$. Trivialising $TM_{|Y}$ yields thus a map from $\mc{S}^2$ to the set $\mc{P}_+$ of $\psi$--positive planes in $\R^7$.

\begin{prp}
The set of $\psi$--positive planes $\mc{P}_+\subset G_4(\R^7)$ is homotopy equivalent to $G_{\psi_0}(\R^7)\cong G_2/SO(4)$, the set of coassociatives. In particular, $\pi_2(\mc{P}_+)\cong\Z_2$.   
\end{prp}
\begin{prf}
Instead of $\mc{P}_+$ we shall consider the dual set $\mc{P}_+^{\perp}\subset G_3(\R^7)$. Restricting $\varphi_0$ to $E\in G_3(\R^7)$ yields a multiple of the volume form induced by $g_{0|E}$, and thus a real number (recall that we always consider the grassmannian of {\em oriented} subplanes). Since $\varphi_0$ is a calibration, we may regard $\varphi_0$ as a map $G_3(\R^7)\to[-1,1]$. By convention we shall orient the orthogonal complements $F^{\perp}$, $F\in\mc{P}_+$, in such a way that $\varphi_0(F^{\perp})<0$, whence $\mc{P}_+^{\perp}=\varphi_0^{-1}\big([-1,0)\big)$. Any fibre $\varphi^{-1}(\cos\alpha)$ is acted on transitively by $G_2$ and contains an element of the form 
$$
E_{\alpha}=e_1\wedge e_2\wedge(\cos\alpha\,e_3+\sin\alpha\,e_4)
$$
with respect to the fixed $G_2$ frame $e_1,\ldots, e_7$ of Section~\ref{groupg2}. To see this, write $E_{\alpha}=x\wedge y\wedge z$ for some unit vectors $x,\,y,\,z\in\R^7$. Since $G_2$ acts transitively on ordered orthonormal pairs with stabiliser $SU(2)$~\cite{hala82}, we may assume that $E_{\alpha}=e_1\wedge e_2\wedge z$ upon transformation by a suitable element in $G_2$. The $SU(2)$--action induced by the inclusion into $G_2$ gives rise to a decomposition $\R^7=\Im\Hh\oplus\Hh$, where $SU(2)$ acts trivially on $\Im\Hh$ and $\Hh=\C^2$ becomes the standard vector representation. We are still free to modify $E_{\alpha}$ without changing $e_1$ and $e_2$ by an element in $SU(2)$. Since this group acts transitively on the unit sphere in $\C^2$, we may transform the unit vector $z=\sum_{i=3}^7z^ie_i$ into $\cos\alpha e_3+\sin\alpha e_4$ with $\cos\alpha=z^3=\varphi_0(E_{\alpha})$. From this one easily deduces that (a) $\pm1$ are the only critical values and (b) that the hessian of $\varphi_0$ is non--degenerate in directions transverse to the orbits $\varphi_0^{-1}(\pm1)$ which are diffeomorphic to $G_2/SO(4)$~\cite{hala82}. Consequently, $\varphi_0$ defines a $G_2$--invariant Morse function in the sense of~\cite{ma89}. By a theorem in the same paper, one can conclude -- in analogy with classical Morse theory -- that $\varphi_0^{-1}\big((-\infty,0]\big)$ is (equivariantly) homotopy equivalent to the disk bundle $G_2\times_{SO(4)} D^4$ (inside the normal bundle) over $\varphi^{-1}\big((-\infty,-1]\big)\cong G_2/SO(4)$ attached. Hence $\mc{P}_+^{\perp}$ is homotopy equivalent to the open disk bundle which can be retracted to the base. In particular, $\mc{P}_+$ is of the same homotopy type as $G_2/SO(4)$. Since $\pi_k(G_2)=0$ for $k=1,\,2$ and $\pi_1\big(SO(4)\big)=\Z_2$, the exact homotopy sequence for fibrations yields the asserted homotopy group.
\end{prf}

\begin{definition}
{\rm Let $D$ be an embedded associative $3$--disk in some topological $G_2$--manifold. We refer to the natural class of $TX_{|\mc{S}^2}$ in $\pi_2(\mc{P}_+)\cong\Z_2$ as the $G_2$--{\em Maslov index} of $D$ and denote it by $\mu_{G_2}(X,D)$.}
\end{definition}

\begin{prp}
If $D\subset M$ is an embedded associative $3$--disk with boundary in some coassociative $X$, then 
$$
\mu_{G_2}(X,D)=\int_{\mc{S}^2}c_1(\nu_X)\mod2=\big(\ind(X,D)+1\big)\mod2.
$$
\end{prp}
\begin{prf}
The natural complex structure of $\mc{S}^2$ gives a natural identification with $\CP^1$ as a complex manifold. Consequently, any complex line bundle is of the form $\mc{O}(n)$ where $n=\int_{\CP^1}c_1\big(\mc{O}(n)\big)$. For instance $T\mc{S}^2=\mc{O}(2)$, and in particular, we find $TX_{|\mc{S}^2}=\mc{O}(2)\oplus\mc{O}(n)$. Let $G_{\psi_0}(\R^7)$ denote the set of coassociative planes which is diffeomorphic to $G_2/SO(4)$. By assumption, there is a map
$$
f:p\in\mc{S}^2\mapsto F_p\in G_{\psi_0}(\R^7).
$$ 
The homotopy class $[f]$ generates $\pi_2\big(G_2/SO(4)\big)$ if and only if its boundary $\partial[f]$ generates $\pi_1\big(SO(4)\big)\cong\Z_2$. If $K:(D^2,\mc{S}^1)\to(\mc{S}^2,N)$ denotes the collapsing map sending $S^1$ to the north pole $N$ of $\mc{S}^2$, a representative of $\partial[f]$ is obtained by restricting the lift $\widetilde{f}:(D^2,S^1)\to G_2$ of the composition $f\circ K$ to $\mc{S}^1$. We can think of $\widetilde{f}_{|\mc{S}^1}\in SO(4)$ as taking values in the set of bases of $F_N$. Its action on $F_N$ can thus be homotoped into the action of the maximal torus $T^2=\mc{S}^1\times\mc{S}^1$ of $SO(4)$ (seen as structure group of $TX_{|\mc{S}^2}$), namely
$$
t\in\mc{S}^1\mapsto\left(\begin{array}{cc} e^{2it} & 0 \\ 0 & e^{nit}\end{array}\right).
$$
On the other hand, the inclusion of the torus induces an epimorphism $\pi_1(T^2)=\Z\oplus\Z\to\pi_1\big(SO(4)\big)=\Z_2$ which is $\mod 2$ reduction of the second summand, whence $[\partial f]=n\mod 2$.
\end{prf}

\begin{ex}
We consider again the coassociative germ example of Section~\ref{indexa}, where $Y$ is the $3$--disk $D$ inside $\Im\Hh\subset\Im\Oc$. Then $\ind(X,D)=1$. On the other hand, $TX_{|S^2}=\mc{O}(2)\oplus\mc{O}(0)$ from which we conclude $\mu_{G_2}(X,D)=0$ in accordance with the proposition. 
\end{ex}

\newpage

\begin{center}
{\bf Acknowledgements}
\end{center}
The first author was supported by the French Agence nationale pour la recherche whilst working on this project. He would like to thank Vincent Borrelli, Jean--Claude Sikorav and Jean--Yves Welschinger for their interest and comments on this work. The second author was partially supported by the SFB 647 ``Space.Time.Matter'' funded by the DFG. He thanks Bernd Ammann, Uli Bunke and Dominic Joyce for useful discussions and suggestions.

\bigskip

\noindent 
D.~\textsc{Gayet}:
Institut Camille Jordan, CNRS UMR 5208, UFR de math\'ematiques, Universit\'e Lyon I B\^atiment  Braconnier, F--69622 Villeurbanne Cedex, France\\
e-mail: \texttt{gayet@math.univ-lyon1.fr}

\medskip

\noindent
F.~\textsc{Witt}: Mathematisches Institut der Universit\"at M\"unster, Einsteinstra{\ss}e 62, D--48149 M\"unster, F.R.G.\\
e-mail: \texttt{frederik.witt@uni-muenster.de}
\end{document}